# Groups of ribbon knots

Ka Yi Ng[*]

February 8, 1995
Revised July 7, 1995

**Abstract**

We prove that for each positive integer $n$, the $V_n$-equivalence classes of ribbon knot types form a subgroup $\mathcal{R}_n$, of index two, of the free abelian group $\mathcal{V}_n$ constructed by the author and Stanford. As a corollary, any non-ribbon knot whose Arf invariant is trivial cannot be distinguished from ribbon knots by finitely many independent Vassiliev invariants. Furthermore, all non-trivial additive knot cobordism invariants are not of finite type. We prove a few more consequences about the relationship between knot cobordism and $V_n$-equivalence of knots. As a bi-product, we prove that the number of independent Vassiliev invariants of order $n$ is bounded above by $\frac{(n-2)!}{2}$ if $n > 5$, improving the previously known upper bound of $(n-1)!$.

## Introduction

In 1990, Vassiliev [V] introduced the space of finite type invariants of knots which arose from viewing knots as embeddings of $S^1$ into $\mathbf{R}^3$. This broadens the classical approach of studying knots. Subsequent authors ([BL], [B-N], [Ko], [L2], [L3]) further studied Vassiliev's invariants and made connections between these invariants and the existing knot invariants. It was shown ([BL], [Tr], [D]) that some knot invariants such as the unknotting number, the signature, the genus and the crossing number of a knot are not finite type invariants. Whether Vassiliev's invariants distinguish all knots is still an open question.

[*]The author is a graduate student at Columbia University, working under the supervision of Joan Birman, and this work is part of her Ph.D. thesis. Partial support by NSF Grant DMS-94-02988 is acknowledged.

In 1994, Gusarov [G] defined an $n$-equivalence relation on knots and showed that the $n$-equivalence classes of knot types determine an abelian group $\mathcal{G}_n$ for each positive integer $n$, where the group operation is induced by the operation of taking connected sums. Inspired by his work and based on Vassiliev's work, the author and Stanford introduced the notion of '$V_n$-equivalence' of knots and proved that such relation also gives a group $\mathcal{V}_n$ of knots for each positive integer $n$. The group $\mathcal{V}_n$ is free abelian and its rank is the number of linearly independent primitive Vassiliev invariants [B-N] of order $\leq n$. The quotient group of $\mathcal{G}_n$ by its torsion subgroup is isomorphic to the group $\mathcal{V}_n$. The relationship between cobordism, $n$-equivalence and $V_n$-equivalence of knots were also discussed in [NS].

In this paper, we continue to study the relationship between cobordism and $V_n$-equivalence of knots. Motivated by the question whether knot cobordism invariants are finite type invariants, we look into the ribbon knots and show that the semigroup of ribbon knots can be made into a group by using the $V_n$-equivalence relation. To state the main result precisely, we first give a few definitions. Two knots $K$, $K'$ are said to be $V_n$-*similar* if $v(K) = v(K')$ for any primitive rational Vassiliev invariant $v$ of order $\leq n$. One can easily check that the $V_n$-equivalence is an equivalence relation and the $V_n$-equivalence classes of knots form a semigroup $\mathcal{V}_n$, the operation $\#$ being induced by the operation of taking connected sums. In [NS], it was shown that $\mathcal{V}_n$ forms a group. Define a knot $K'$ to be $V_n$-*inverse* to a knot $K$ if $v(K') = -v(K)$ for any primitive rational Vassiliev invariant $v$ of order $\leq n$. Let $\mathcal{R}_n$ be the semigroup of $V_n$-equivalence classes of ribbon knots under the operation $\#$. Then $\mathcal{R}_n$ is a subsemigroup of $\mathcal{V}_n$. The main result in this paper is

**Theorem 3.2** *$\mathcal{R}_n$ forms a subgroup of the free abelian group $\mathcal{V}_n$ of index two. So its rank is the same as the rank of $\mathcal{V}_n$ and is the number of linearly independent primitive invariants of order $\leq n$.*

The proof of Theorem 3.2 follows the approach described in [NS]. We extend the sequence of ribbon knots described in [Kan] and [O1] to a sequence of families $\{K_\sigma\}_{\sigma \in O_n}$ of ribbon knots and make use of the extended classes of ribbon knots to construct $V_n$-inverses of all ribbon knots.

Theorem 3.2 has several interesting consequences about the relationship between knot cobordism and the $V_n$-equivalence relation. An immediate consequence is that finitely many independent Vassiliev invariants are not sufficient to distinguish the ribbon knots from non-ribbon knots (Corol-lary 3.3). This does not answer the open question of whether Vassiliev invariants distinguish all knots, but this already provides some new and interesting information.

Another consequence of Theorem 3.2 is that all non-trivial additive rational knot cobordism invariants are not of finite type (Corollary 3.4). In [L2], Lin showed that the $\mu$-invariants, introduced by Milnor ([M]), are all of finite type. These $\mu$-invariants are invariants of concordance classes of string links ([L1]). Therefore, the cobordism relations on knots and links are quite different.

The following corollary gives a connection between knot cobordism and $V_n$-equivalence of knots.

**Corollary 3.5** *Let $K$ and $K'$ be two knots with the same Arf invariants. Then for any positive integer $n$, there exists a knot $K_n$ which is cobordant to $K$ and is $V_n$-similar to $K'$.*

In [Tr], Trapp showed that the signature of a knot is not a Vassiliev invariant. This also follows from Corollary 3.4. Thus, the signature is independent of Vassiliev invariants of knots. Indeed, we can give examples of knots with the same Vassiliev invariants up to an arbitrarily high order but with an arbitrary signature (Corollary 3.6).

In [NS], the $V_n$-inverses of knots constructed by pure braid commutators are also the $n$-inverses in the sense of Gusarov. By using the ribbon knots $\{K_\sigma\}_{\sigma \in O_n}$, we show that the $n$-equivalence classes of ribbon knots give a subgroup of the group $\mathcal{G}_n$ for each integer $n$ (Theorem 4.1). Thus, we have a stronger result that any non-trivial additive knot cobordism invariant which takes values in an arbitrary group cannot be of finite type. Furthermore, the Arf invariant is the only Vassiliev invariant which is also an additive knot cobordism invariant.

Along the proof of Theorem 3.2, we are able to give an upper bound for the number of linearly independent primitive Vassiliev invariants of order $n$ (Theorem 5.2). By using these estimates, we compute that the dimension of the space $V_n$ of all Vassiliev invariants of order $\leq n$ modulo the space $V_{n-1}$ of invariants of order $< n$ is bounded above by $\frac{(n-2)!}{2}$ if $n > 5$ (Corollary 5.3). They are finer than the upper bounds given by Chmutov and Duzhin ([CD1]).

This paper is organized as follows. Section 1 gives a brief introduction to Vassiliev invariants and proves Theorem 1.3 about $n$-gons, which will be used in proving Theorem 3.2 and Theorem 5.2. Section 2 describes the family $\{K_\sigma\}_{\sigma \in Y_n}$ of ribbon knots and computes their Vassiliev invariants of order up to $n$. In Section 3, we prove Theorem 3.2 and its consequences. We prove Theorem 4.1 in Section 4. Finally, Section 5 gives an upper bound for the dimension of the quotient space $V_n/V_{n-1}$ for $n > 5$.

## Acknowledgement

The author is deeply indebted to her advisor, Joan Birman, and thanks for her comments, suggestions and for suggesting the rephrased statement (Corollary 3.3) of Theorem 3.2. Thanks also goes to T. Stanford, who helps a lot in understanding the subject. The author wishes to thank her husband, T. S. Fung, for his support and helpful discussions. Thanks also goes to H. Bass, M. Greenwood and X. S. Lin for many discussions. Thanks goes to Thang Le for reading this manuscript in great details and for his helpful suggestions and advice.

## 1 Chord diagrams and $n$-gons

For completeness, we include here the definition of Vassiliev invariants, following Birman and Lin ([BL]).

*Singular knots* are immersions of $S^1$ in $\boldsymbol{R}^3$ whose singularities are transverse double points, considered up to the equivalence relation generated by flat vertex isotopy ([BL]).

A *Vassiliev invariant $v$ of order $n$* is a knot invariant which takes values in some abelian group and which can be extended to an invariant of singular knots by the relation depicted in Figure 1 such that as a singular knot invariant, $v(K_D) = 0$ for any singular knot $K_D$ with more than $n$ double points and $v(K_C) \neq 0$ for some singular knot $K_C$ with $n$ double points. A rational Vassiliev invariant $v$ is one which takes rational values and if $v$ is additive under connected sums of knots, $v$ is said to be *primitive*. In this paper, we shall only consider invariants which vanish on the unknot.

A *chord diagram of order $n$* is a counterclockwise oriented circle with $n$ chords whose endpoints are all distinct ([BL], [B-N]). A *split* chord diagram is a chord diagram whose set of chords can be separated by a dotted line into two non-empty disjoint subsets. A chord diagram is *non-split* if it is not split.

To each singular knot with $n$ double points, we may associate a chord diagram of order $n$ by connecting the preimages of each double point of the immersion with a chord. Given a Vassiliev invariant $v$ of order $n$, since $v(K_D) = 0$ for any singular knot $K_D$ with $(n+1)$ double points, $v$ is well-defined on chord diagrams of order $n$. Indeed, it can be shown that $v$ is determined by a set of axioms and initial data which is in the form of a table, called an actuality table, containing values for the chord diagrams and immersions that respect them. Modulo the invariants of order $< n$, invariants of order $n$ are determined by their values on chord diagrams of order $n$. (For details, see [BL].)

In [B-N], Bar-Natan generalised the concept of chord diagrams and introduced the Chinese Character Diagrams. A *Chinese Character Diagram* (CCD) of order $n$ is a connected graph of $2n$ vertices which consists of a distinguished subgraph, homeomorphic to a circle. The circle is called the *external circle* and is oriented while other edges are all unoriented. Each vertex is trivalent and has an orientation, that is, a cyclic ordering of the edges incident to it. Vertices on the external circle are called *external* and the others are called *internal*. We call the subgraph obtained by deleting the external circle from the CCD the *internal graph*. CCD's are considered up to isomorphisms that respect the above structures. They are usually drawn as in Figure 2, with the convention that the external circle and the vertices are oriented counterclockwise.

Let $B_n$ be the additive group generated by the CCD's of order $n$ modulo the STU relations shown in Figure 3 and let $B'_n$ be its 2-torsion subgroup. Let $A_n$ be the additive group generated by chord diagrams of order $n$ modulo the 4T relations shown in Figure 4. Then the inclusion of chord diagrams into CCD's induces an isomorphism between the abelian groups $A_n$ and $B_n$ ([B-N]). The STU relations generate the IHX relations, shown in Figure 5, and the antisymmetry of internal vertices, shown in Figure 6.

Several authors, for example, Bar-Natan, Chmutov-Duzhin, Greenwood, Melvin, Piunhikan have noticed that the addition of a chord of length two (that is, adding a chord such that the newly added chord has only one endpoint between its endpoints) to a non-split chord diagram modulo the split diagrams is independent of where the chord of length two is placed.

**Theorem 1.1** *Let $C_n$ be the subgroup of $B_n$ spanned by connected CCDs. Then there exists a group homomorphism $\phi_n : C_n \to B_{n+1}/B'_{n+1}$ which is defined by adding a chord of length two.*

We define an *incomplete $n$-gon* to be a CCD with $n$ external vertices and $n$ internal vertices, whose internal graph consists of one and only one cycle (without any backtracking), and the cycle is of length $k < n$. Figure 2 shows an example of an incomplete 5-gon. If $k = n$, then we call the diagram a *complete $n$-gon*.

In this section, we show

**Theorem 1.2** *For $n > 2$, $B_n$ is generated by complete $n$-gons, split diagrams and its 2-torsion subgroup $B'_n$.*

In [NS], a *one-branch tree diagram of order $n$* was defined as a CCD whose internal graph is a tree isomorphic to the standard $n$-tree shown in Figure 7 preserving the vertex orientations. Such a one-branch tree diagram $T$ of order $n$ determines a permutation $\sigma \in S_n$ as follows. Label the branches of the standard tree as in Figure 7. Under the isomorphism between the standard tree and the internal graph $G$ of $T$, the branches of $G$ are labelled. Number the $(n+1)$ external vertices of $T$ by $0, 1, 2, \ldots, n$ in counterclockwise direction such that the end of branch 0 of $G$ is numbered by 0. Then the correspondence between the branches of $G$ and their endpoints determines a permutation $\sigma \in S_n$. Namely, branch $i$ of $G$ has its end labelled as $\sigma(i)$ for $i = 1, 2, \ldots, n$. Conversely, we may construct a unique one-branch tree diagram of order $n$ from a permutation $\sigma \in S_n$ in the same way. We shall denote the diagram by $T_\sigma$.

In [NS], it was shown that

**Lemma 1.3** [NS] *$B_n$ is generated by split diagrams and one-branch tree diagrams of order $n$.*

Therefore, to show Theorem 1.2, we only need the following lemma.

**Lemma 1.4** *Any one-branch tree diagram in $B_n$ ($n > 2$) can be expressed as an integral sum of complete $n$-gons and diagrams in $B'_n$.*

**Proof.** Let $T_\sigma$ be the one-branch tree diagram determined by the permutation $\sigma \in S_n$ for $n > 2$. We first claim that $T_\sigma$ in $B_n$ is equal to an integral sum of complete or incomplete $n$-gons, together with possibly some diagrams in $B'_n$. We shall prove this claim by induction on $k = |\sigma(n-1) - \sigma(n-2)|$.

Step 1:

Let $k = |\sigma(n-1) - \sigma(n-2)| = 1$. Without loss of generality, assume that $\sigma(n-1) = \sigma(n-2) + 1$. Let $\gamma_\sigma$ be the oriented path along the circle from vertex $\sigma(n)$ to vertex $\sigma(n-2)$ and let $l(\gamma_\sigma)$ be its length.

If $l(\gamma_\sigma) = 1$, then we may apply the IHX relation to the portion of the graph containing the external vertices $\sigma(n-2)$, $\sigma(n-1)$ and $\sigma(n)$ as shown in Figure 8. If we apply the STU relation to the vertex $v$ in the middle diagram and to the vertex $v'$ in the right diagram of Figure 8, we find that two of the diagrams so-obtained are split, and the other two differ by the placement of a chord of length two. Theorem 1.1 implies that $T_\sigma$ is trivial in $B_n/B'_n$.

If $l(\gamma_\sigma) = r > 1$, then let $\sigma(k_0), \sigma(k_1), \ldots, \sigma(k_r)$ be the labels of the external vertices along the oriented path $\gamma_\sigma$. Apply the STU relation to the two edges labelled as $k_0(=n)$ and $k_1$ as shown in Figure 9. Then $T_\sigma$ is a sum of an (complete or incomplete) $n$-gon and a one-branch tree diagram $T_{\sigma_1}$ of order $n$, where $\sigma_1(n) = k_1, \sigma_1(k_1) = n$ and $\sigma_1(k) = \sigma(k)$ for $k \neq n, k_1$. Then $l(\gamma_{\sigma_1}) = r - 1$. Repeat the application of the STU relation in a similar way. We will finally get a signed sum of $(r-1)$ complete or incomplete $n$-gons and a one-branch tree diagram $T_{\sigma'}$ with $l(\gamma_{\sigma'}) = 1$. As we have discussed, $T_{\sigma'} = 0$ in $B_n/B'_n$. So the claim for $k = 1$ is proved.

Step 2:

Suppose $k = |\sigma(n-1) - \sigma(n-2)| > 1$. Again we may assume that $\sigma(n-1) > \sigma(n-2)$. We may further assume that $\sigma(n)$ is not between $\sigma(n-2)$ and $\sigma(n-1)$ because if $\sigma(n-2) < \sigma(n) < \sigma(n-1)$, then up a change of sign of $T_\sigma$, we may replace $\sigma$ by $\sigma_1$ where $\sigma_1(n-1) = \sigma(n), \sigma_1(n) = \sigma(n-1)$ and $\sigma_1 = \sigma$ on $\{1, 2, \ldots, n-2\}$. Now apply the STU relation to the branch $(n-1)$ and the one between branches $(n-1)$ and $(n-2)$, right next to branch $(n-1)$. As shown in Figure 10, $T_\sigma$ is expressed as a sum of an (complete or incomplete) $n$-gon and a one-branch tree diagram $T_{\sigma'}$ where $|\sigma'(n-1) - \sigma'(n-2)| = k-1$. By the induction hypothesis, $T_{\sigma'}$ is an integral sum of complete or incomplete $n$-gons, and some diagrams in $B'_n$. The claim is proved.

It remains to show that any incomplete $n$-gon is an integral sum of complete $n$-gons. Let $D$ be an incomplete $n$-gon. Then $D$ contains a cycle of length $k < n$ joining some internal vertices. Pick an internal vertex $v_1$ in the cycle which is connected by an edge to another internal vertex $v_2$ not in the cycle. Apply the IHX relation such that the vertical piece of I corresponds to the edge joining $v_1$ and $v_2$. Then $D$ is expressed as a difference of two $n$-gons both of which contain a cycle of length $(k+1)$. So inductively, $D$ is generated by complete $n$-gons. ∥

Combining Lemma 1.3 and Lemma 1.4, Theorem 1.2 follows immediately.

For $n > 1$, define a map $f$ from $S_n$ to the set of all complete $n$-gons as follows. Let $G_n$ be the graph which is a regular $n$-gon with an edge attached to each of its $n$ vertices. Label the edges by $1, \ldots, n$ in the counterclockwise direction. The graph $G_5$ is shown in Figure 11. To each $\sigma \in S_n$, we assign a complete $n$-gon whose internal graph $G$ is isomorphic to the graph $G_n$, preserving the vertex orientations, and if the external vertices are numbered as $1, \ldots, n$ in the counterclockwise direction, then the edge $i$ of $G$ should be attached to the external vertex $\sigma(i)$ for $1 \leq i \leq n$. Call this complete $n$-gon $f(\sigma)$. It is obvious that the map $f$ is well-defined. We shall denote the image of $f$ by $O_n$. By the antisymmetry relation, any complete $n$-gon which is not in $O_n$ is equal to $\pm f(\sigma)$ for some $\sigma \in S_n$. Then by Theorem 1.2, $B_n$ is generated by split diagrams, complete $n$-gons in $O_n$ and its 2-torsion subgroup.

We remark that the map $f : S_n \to O_n$ is not one-to-one. Define an ordering $\leq$ on $S_n$ such that for any $\sigma$, $\sigma'$ in $S_n$, $\sigma \leq \sigma'$ if $\sigma = \sigma'$ or for some integer $k$ with $1 \leq k \leq n$, $\sigma(k) < \sigma'(k)$ and $\sigma(j) = \sigma'(j)$ for $1 \leq j < k$. If among all the preimages of a given complete $n$-gon $D \in O_n$, we call the smallest element with respect to the above ordering $\leq$ a *canonical representative of $D$*, then every complete $n$-gon in $O_n$ has a unique canonical representative.

## 2 The Ribbon knots $K_\sigma$

A *ribbon knot* in $\mathbf{R}^3$ is a knot which bounds an immersed disk in $\mathbf{R}^3$ with ribbon intersections only. That is, it is the boundary of the image of an immersion $j$ of a disk $D$ in $\mathbf{R}^3$ where the intersections of the immersed disk are all transverse and the double point set of $j$ consists of finitely many pairs of mutually disjoint arcs $(I_\alpha, I'_\alpha)$ in $D$ such that $j(I_\alpha) = j(I'_\alpha)$ and $I_\alpha$ is properly embedded in $D$ and $I'_\alpha$ is contained in $\text{int}(D)$.

In this section we shall describe the ribbon knots $K_\sigma$ associated to the complete $n$-gons $f(\sigma)$ in $O_n$ such that the CCD $f(\sigma)$ and the knot $K_\sigma$ have the same primitive Vassiliev invariants of order $n$ and $K_\sigma$ has trivial Vassiliev invariants of order less than $n$.

Let $\sigma$ be the canonical representative of a complete $n$-gon in $O_n$. Then $\sigma(1)$ must be one. We shall construct a ribbon knot $K_\sigma$ by drawing a projection of $K_\sigma$ on $\mathbf{R}^2$ as follows.

Let $b_i = \sigma^{-1}(i)$ for $1 \leq i \leq n$. Consider a dotted circle with $2n$ marked points labelled as $x_{11}, x_{12}, \ldots, x_{n,1}, x_{n,2}$ in the counterclockwise direction (Figure 12). Starting from the point $x_{n,2}$, draw lines from $x_{b_i-1,2}$ to $x_{b_{i+1}-1,1}$ in ascending order of $i$ for $i = 1,\ldots,n$ (where $x_{0i}$ is defined as $x_{ni}$) such that the first passage across a double point in the projection is always an overpass. Call the resulted diagram $P$. Note that if we add the arcs joining $x_{i1}$ and $x_{i2}$ for $1 \leq i \leq n$ to the diagram, then we get a layered diagram of the unknot and it bounds a disk $D$ in $\mathbf{R}^3$. Now put the diagram $P$ over the region $P'$ in Figure 13 such that the points $x_{ij}$'s in the diagram $P$ coincide with the corresponding points $x_{ij}$'s on the boundary of the region $P'$. Then the diagram obtained is a projection of a ribbon knot. Denote the projection by $P_\sigma$ and the ribbon knot by $K_\sigma$. An example is shown in Figure 14 which is the ribbon knot constructed from (34). The immersed disk that the knot bounds is shaded.

The ribbon knots we constructed here are generalisations of the ribbon knots described by Kanenobu. The knots $K_\sigma$ for $\sigma = (1, 2, \ldots, n)$ are precisely the examples $K(0, \ldots, 0)$ given in §5 in [Kan].

Ohyama [O2] and Gusarov [G] introduced the concept of $n$-triviality of a knot. Let $S_1, S_2, \ldots, S_{n+1}$ be $(n+1)$ mutually disjoint sets of crossings chosen from a regular projection $P$ of a knot $K$. We call $\boldsymbol{S} = (S_1, S_2, \ldots, S_{n+1})$ an $(n + 1)$-*scheme in the projection $P$ of the knot $K$*. For each $\boldsymbol{x} = (x_1, \ldots, x_{n+1}) \in \{0, 1\}^{n+1}$, denote by $K(\boldsymbol{x}, \boldsymbol{S})$ the knot obtained by changing the crossings in $S_i$'s for which $x_i = 1$. Following Gusarov [G], $K$ is said to be $n$-*trivial* if there exists an $(n + 1)$-scheme $\boldsymbol{S}$ in a projection of $K$ such that $K(\boldsymbol{x}, \boldsymbol{S})$ is isotopic to the unknot for any non-zero $\boldsymbol{x}$ in $\{0, 1\}^{n+1}$.

In the projection $P_\sigma$ of the ribbon knot $K_\sigma$, if we let $T_j$ be the set containing the two crossings $c_{j1}, c_{j2}$ for $1 \leq j \leq n$ and let $\boldsymbol{T}_\sigma = (T_1, \ldots, T_n)$, then $K_\sigma(\boldsymbol{x}, \boldsymbol{T}_\sigma)$ is isotopic to the unknot for any non-zero $\mathbf{x} \in \{0, 1\}^n$. This shows that the ribbon knots $K_\sigma$'s are $(n-1)$-trivial.

By using the formula given by Ohyama in [O2], we may compute the Vassiliev invariants of the ribbon knot $K_\sigma$.

**Proposition 2.1** *Let $\sigma \in O_n$ be a canonical representative where $n \geq 2$. Then $K_\sigma$ has trivial Vassiliev invariants of order less than $n$. For any Vassiliev invariant $v$ of order $n$, $v(K_\sigma) = v(f(\sigma))$.*

**Proof.** Apply Lemma 3 in [O2] to the ribbon knot $K_\sigma$ equipped with the scheme $\boldsymbol{T}_\sigma$. Then for any Vassiliev invariant $v$,

$$v(K_\sigma) = \sum_{\substack{i_j = 1,2 \\ 1 \leq j \leq n}} \varepsilon_{i_1}^{(1)} \cdots \varepsilon_{i_n}^{(n)} v(K_\sigma(i_1, \ldots, i_n))$$

where $\varepsilon_1^{(j)} = -\varepsilon_2^{(j)} = 1$ for $1 \leq j \leq n$, and $K_\sigma(i_1, \ldots, i_n)$ denotes the singular knot obtained from $K_\sigma$ by changing the first $(i_j - 1)$ crossings of $S_j$ and smashing the $i_j$th crossing into a double point for $1 \leq j \leq n$.

The first statement of the proposition follows immediately because the evaluation of invariants of order $< n$ on the singular knots $K_\sigma(i_1, \ldots, i_n)$ are zero. The second part of the proposition will be shown as follows.

Consider a counterclockwise oriented circle marked with $3n$ points labelled as $z_1^{(b_1-1)}, y_{b_1}, z_2^{(b_1-1)}, z_1^{(b_2-1)}, y_{b_2}, z_2^{(b_2-1)}, \ldots, z_1^{(b_n-1)}, y_{b_n}, z_2^{(b_n-1)}$ (Figure 15). For each $(i_1, \ldots, i_n) \in \{1,2\}^n$, join the points $y_j$ and $z_{i_j}^{(j)}$ with a chord $L_{i_j}^{(j)}$ for $1 \leq j \leq n$. Call the chord diagram $C(i_1, \ldots, i_n)$. Then $C(i_1, \ldots, i_n)$ is the chord diagram that the singular knot $K_\sigma(i_1, \ldots, i_n)$ respects. Therefore, to complete the proof, it suffices to show that in $B_n$,

$$\sum_{\substack{i_j=1,2 \\ 1 \leq j \leq n}} \varepsilon_{i_1}^{(1)} \cdots \varepsilon_{i_n}^{(n)} C(i_1, \ldots, i_n) = f(\sigma) \tag{1}$$

Note that if we omit the chord $L_{i_1}^{(1)}$ from the diagram $C(i_1, \ldots, i_n)$ for $i_1 = 1, 2$, then the two chord diagrams of order $n-1$ are the same. Thus, if we apply the STU relation to the arc containing the points $z_1^{(1)}, y_2, z_2^{(1)}$ as shown in Figure 16, the difference of the two diagrams $C(1, i_2, \ldots, i_n)$ and $C(2, i_2, \ldots, i_n)$ is a CCD which we call $C(*, i_2, \ldots, i_n)$. Then $C(*, i_2, \ldots, i_n)$ is the diagram whose internal graph consists of the chords $L_{i_j}^{(j)}$'s joining the points $y_j$ and $z_{i_j}^{(j)}$ for $2 < j \leq n$ and a tree of order 2 isomorphic to the standard 2-tree (preserving the vertex orientations), with branches $0, 1, 2$ attached to the points $y_1, y_2$ and $z_{i_2}^{(2)}$ respectively. If $n = 2$, we apply the STU relation to the diagrams $C(*, 1)$ and $C(*, 2)$ around the arc containing the points $z_1^{(2)}, y_1, z_2^{(2)}$ to obtain a complete 2-gon. Otherwise, we apply the STU relation to the two diagrams $C(*, 1, i_3, \ldots, i_n)$ and $C(*, 2, i_3, \ldots, i_n)$ around the arc containing the points $z_1^{(2)}, y_3, z_2^{(2)}$. The difference of these two diagrams is the CCD whose internal graph consists of the chords $L_{i_j}^{(j)}$'s joining the points $y_j$ and $z_{i_j}^{(j)}$ for $3 < j \leq n$ and a tree of order 3 isomorphic to the standard 3-tree (preserving the vertex orientations), with branches $0, 1, 2, 3$ attached to the points $y_1, y_2, y_3$ and $z_{i_3}^{(3)}$ respectively. Repeat the procedure until the signed sum on the left hand side of (1) is combined to one CCD. Then the external vertices of the CCD are labelled as $y_{b_1}, y_{b_2}, \ldots, y_{b_n}$ and the internal graph $G$ is isomorphic to the graph $G_n$, preserving the vertex orientations. The branch $i$ of $G$ is attached to the external vertex $y_i$ for $1 \leq i \leq n$. By relabelling each external vertex $y_{b_i}$ as $i$, the branch $i$ of $G$ becomes attached to the vertex $\sigma(i)$. Therefore, it is identical to $f(\sigma)$. ∥

We may also associate a ribbon knot $K_\sigma^{-1}$ to $-f(\sigma)$.

**Proposition 2.2** *For any canonical representative $\sigma \in O_n$ ($n \geq 2$), there exists a ribbon knot $K_\sigma^{-1}$ such that $v(K_\sigma^{-1}) = -v(K_\sigma)$ for any Vassiliev invariant $v$ of order $\leq n$.*

**Proof.** In the projection $P_\sigma$ of the ribbon knot $K_\sigma$, change the local picture containing crossings $c_{11}$ and $c_{12}$ as shown in Figure 17. Let $T_1'$ be the set containing the two crossings $c_{11}'$ and $c_{12}'$. By computing the Vassiliev invariants of this new ribbon knot with the scheme $(T_1', T_2, \ldots, T_n)$, the proposition follows easily. ∥

## 3 The group $\mathcal{R}_n$ of ribbon knots

In this section, we prove Theorem 3.2 and its consequences. Unless otherwise stated, the Vassiliev invariants mentioned here are rational. Recall that $\mathcal{R}_n$ is the abelian semigroup of $V_n$-equivalence classes of ribbon knots with the semigroup operation $\#$.

**Lemma 3.1** *Let $D_1, \ldots, D_r$ be all the distinct non-split chord diagrams of order $n > 2$. Let $\lambda_1, \ldots, \lambda_r$ be $r$ integers. Then there exists a ribbon knot $K(\lambda_1, \ldots, \lambda_r)$ whose invariants of order $< n$ are zero and $v(K) = \sum_{i=1}^{r} \lambda_i v(D_i)$ for any primitive rational Vassiliev invariant $v$ of order $n$.*

**Proof.** Let $v$ be a primitive Vassiliev invariant of order $n$. Since the evaluation of $v$ on any split diagram is zero [NS], by Theorem 1.2 and the remark after the proof of Theorem 1.4, $\sum_{i=1}^{r} \lambda_i v(D_i) = \sum_{\sigma \in I} \epsilon_\sigma \mu_\sigma v(f(\sigma))$ for some positive integers $\mu_\sigma$ for $\sigma \in I \subseteq O_n$ and $\epsilon_\sigma \in \{\pm 1\}$.

Given any knot $K$, denote by $mK$ the connected sum of $m$ copies of the knot $K$. It follows from Proposition 2.1 and Proposition 2.2 that the knot $K = \#_{\sigma \in I} \mu_\sigma K_\sigma^{\epsilon_\sigma}$ and the integral sum $\sum_{\sigma \in I} \epsilon_\sigma \mu_\sigma f(\sigma)$ of complete $n$-gons have the same primitive rational Vassiliev invariants of order $n$. The knot $K$ has trivial Vassiliev invariants of order $< n$ since the connected sum of two $(n-1)$-trivial knots is $(n-1)$-trivial and any $(n-1)$-trivial knot has trivial Vassiliev invariants of order $< n$ ([G]). Therefore, the lemma is proved.∥

**Theorem 3.2** $\mathcal{R}_n$ *forms a subgroup of the free abelian group* $\mathcal{V}_n$ *of index two. So its rank is the same as the rank of* $\mathcal{V}_n$ *and is the number of linearly independent primitive invariants of order* $\leq n$.

**Proof.** First observe that if $K$ is a knot whose Vassiliev invariants of order $< n$ are trivial, then there exist rational numbers $\lambda_i$ and chord diagrams $D_i$ of order $n$ such that $v(K) = \sum_{i=1}^{r} \lambda_i v(D_i)$ for any Vassiliev invariant $v$ of order $n$. ([BL], [B-N], [V])

In [BL] and [B-N], it is shown that the coefficient $a_2(K)$ in the Conway polynomial of a knot $K$ is a Vassiliev invariant of order 2. In [Kau], Kauffman showed that $a_2(K) \equiv \mathrm{Arf}(K)$ (mod2). Therefore, $a_2$ modulo 2 is a knot cobordism invariant. Since the space of Vassiliev invariants of order 2 is of dimension one, any Vassiliev invariant of order 2 modulo 2 is a knot cobordism invariant.

Let $K$ be a knot. We shall show by induction on $n \geq 2$ that the ribbon knot $K$ has an $V_n$-inverse in $\mathcal{R}_n$. ($\mathcal{R}_1 = \mathcal{V}_1$ is trivial.)

Consider the case $n = 2$. If $v(K) = 0$ for some invariant $v$ of order 2, then $v(K) = 0$ for all invariants $v$ of order 2 and $K$ itself is an $V_n$-inverse of $K$.

Suppose $v(K) \neq 0$ for all invariants $v$ of order 2. Then following the above discussion, there exists an integer $\lambda$ such that $v(K) = 2\lambda v(D)$ for any Vassiliev invariant $v$ of order 2. Note that by the STU relations, the complete 2-gon $f(\sigma)$ is equal to $2D$ in $B_2$ where $\sigma = (12)$. By Proposition 2.1 and Proposition 2.2, $v(K) = v(|\lambda|K_\sigma^\epsilon)$ where $\epsilon = \frac{\lambda}{|\lambda|}$ and $|\lambda|K_\sigma^\epsilon$ is the ribbon knot obtained by taking connected sum of $|\lambda|$ copies of the knot $K_\sigma^\epsilon$. Hence, $|\lambda|K_\sigma^{-\epsilon}$ is an inverse of $K$ in $\mathcal{R}_2$.

Suppose the ribbon knot $K$ has an inverse $K_1$ in $\mathcal{R}_{n-1}$. Then for any primitive Vassiliev invariant $v$ of order $< n$, $v(K \# K_1) = v(K) + v(K_1) = 0$. Using the observation made at the beginning of the proof, there exist rational numbers $\lambda_i$ and chord diagram $D_i$ of order $n$ such that $v(K \# K_1) = \sum_{i=1}^{r} \lambda_i v(D_i)$ for any Vassiliev invariant $v$ of order $n$.

Let $d$ be the common denominator of the $\lambda_i$'s. By Lemma 3.1, there exists a ribbon knot $K(-\lambda_1 d, \ldots, -\lambda_n d)$ whose invariants of order $< n$ are zero and $v(K(-\lambda_1 d, \ldots, -\lambda_n d)) = -\sum_{i=1}^{r} \lambda_i d v(D_i)$ for any primitive Vassiliev invariant $v$ of order $n$. Then $d(K \# K_1) \# K(-\lambda_1 d, \ldots, -\lambda_n d)$ has trivial invariants of order $< n$ and for any primitive invariant $v$ of order $n$,

$$v(d(K \# K_1) \# K(-\lambda_1 d, \ldots, -\lambda_n d)) = dv(K \# K_1) + v(K(-\lambda_1 d, \ldots, -\lambda_n d))$$

$$= d\sum_{i=1}^{r}\lambda_i v(D_i) - \sum_{i=1}^{r}\lambda_i dv(D_i)$$
$$= 0$$

Thus, $(d-1)K\#dK_1\#K(-\lambda_1 d,\ldots,-\lambda_n d)$ is an inverse of $K$ in $\mathcal{R}_n$. It follows from the inductive proof that $\mathcal{R}_n$ is of index two. The second statement of the theorem is immediate. $\|$

Theorem 3.2 can be put in the following way:

**Corollary 3.3** *Let $K$ be a non-ribbon knot whose Arf invariant is trivial. Then there is an infinite sequence $\{K_n\}_{n=2}^{\infty}$ of ribbon knots such that for each $n \geq 2$, the given knot $K$ cannot be distinguished from $K_n$ by rational Vassiliev invariants of order at most $n$.*

We may deduce from Corollary 3.3 that

**Corollary 3.4**

**(i)** *All non-trivial additive rational knot cobordism invariants are not of finite type.*

**(ii)** *Any primitive rational Vassiliev invariant of order $n$ cannot be a non-trivial knot cobordism invariant for any integer $n > 2$.*

**Proof.** For (i), suppose on the contrary that there is a non-trivial additive knot cobordism invariant $\omega$ which is a Vassiliev invariant of order $n$. Let $K$ be a knot such that $\omega(K) \neq 0$. By Corollary 3.3, there is a ribbon knot $K'$ such that $\omega(K') = 2\omega(K)$ or $\omega(K)$ according as $n = 2$ or $n > 2$. Thus, in both cases, $\omega(K') \neq 0$. This contradicts that $\omega$ is a knot cobordism invariant. This proves the first statement.

For (ii), the second statement can be shown by a similar argument. $\|$

It is shown in [NS] that the inclusion map from the set $\mathcal{K}$ of all knot types to the knot cobordism group $\mathcal{C}$ does not factor through the $V_n$-equivalence and the inclusion map from $\mathcal{K}$ to the group $\mathcal{V}_n$ does not factor through the cobordism relation. There seems to be no obvious connection between the $V_n$-equivalence and the cobordism relation. In the next corollary, we try to use Theorem 3.2 to make a bridge between the two equivalence relations.

**Corollary 3.5** *Let $K$ and $K'$ be two knots with the same Arf invariants. Then for any positive integer $n$, there exists a knot $K_n$ which is cobordant to $K$ and is $V_n$-similar to $K'$.*

**Proof.** Let $\epsilon = \text{Arf}(K) = \text{Arf}(K')$. Following the proof of Theorem 3.2, we may construct two ribbon knots $L_n$ and $L'_n$ such that $a_2(L_n) + a_2(K) = a_2(K') - a_2(L'_n) = \epsilon$ and for any primitive Vassiliev invariant $v$ of order $m$ with $2 < m \leq n$, $v(L_n) = -v(K)$ and $v(L'_n) = v(K')$. Then the connected sum $K_n$ of the knots $K$, $L_n$ and $L'_n$ is the desired knot. ∥

The following corollary exhibits a way of how the signature is independent of the Vassiliev invariants.

**Corollary 3.6** *Let $K$ be a knot and let $n$ be a positive integer. Then there exists a sequence $\{K_m\}_{m\in\mathbf{Z}}$ of knots such that for each $m$, the signature of $K_m$ is $2m$ and its Vassiliev invariants of order $\leq n$ match with those of the knot $K$.*

**Proof.** Let $T_m$ be the torus knot of type $(2, 2m+1)$ for any integer $m \geq 1$ and let $T_{-m}$ the mirror image of the knot $T_m$ for $m > 0$. It is known that $\sigma(T_m) = 2m$ for any $m$ and

$$\text{Arf}(T_m) = \begin{cases} 0 & \text{when } |m| \equiv 0, 3 (\text{mod } 4); \\ 1 & \text{when } |m| \equiv 1, 2 (\text{mod } 4). \end{cases}$$

Suppose $\text{Arf}(K) = 0$. Then for $|m| \equiv 0, 3 (\text{mod } 4)$, $\text{Arf}(K) = \text{Arf}(T_m)$. By Corollary 3.5, there exists a knot $K_m$ which is cobordant to $T_m$ and is $V_n$-similar to $K$. Let $K'$ be the figure eight knot. For $|m| \equiv 1, 2 (\text{mod } 4)$, $\text{Arf}(K) = \text{Arf}(T_m \# K')$. So there exists a knot $K_m$ cobordant to $T_m \# K'$ and $V_n$-similar to $K$. A similar method can be applied to the case when $\text{Arf}(K) = 1$. It is obvious that these knots $K_m$'s are the desired ones. ∥

The above arguments show that we may obtain similar results for other additive knot cobordism invariant $\omega$ provided that we have examples of knots with arbitrary values of the invariant $\omega$ and the Arf invariant.

## 4 A subgroup of Gusarov's group $\mathcal{G}_n$

In this section we extend the results in Section 3 and show that the $n$-equivalence classes of ribbon knots form a subgroup of the group $\mathcal{G}_n$ for each integer $n$. Note that two knots are $n$-equivalent if and only if they have the same additive Vassiliev invariants of order $\leq n$ ([G]). Therefore, it is not difficult to see that the $n$-inverses of the $n$-equivalence classes of ribbon knots can be constructed in a similar way as in the proof of Theorem 3.2.

**Theorem 4.1** *The $n$-equivalence classes of ribbon knots form a subgroup $\mathcal{RG}_n$ of the group $\mathcal{G}_n$ of index two.*

**Proof.** In view of Theorem 3.2, we only need to complete the case in the inductive step where the group $B_n$ has 2-torsion.

Suppose a ribbon knot $K$ has an inverse $K_1$ in $\mathcal{RG}_{n-1}$ and $B_n$ has 2-torsion. Then $v(K\#K_1) = 0$ for any additive Vassiliev invariant $v$ of order $< n$. This implies that there exists an integral sum $s$ of chord diagrams of order $n$ such that $K\#K_1$ and this sum $s$ have the same additive Vassiliev invariants of order $n$. By Theorem 1.2, the sum $s \in B_n$ is equal to an integral sum $t$ of complete $n$-gons modulo split diagrams and diagrams in $B'_n$. Apply Proposition 2.1 and Proposition 2.2 to get a ribbon knot $K_t$ whose additive invariants of order $< n$ are trivial and $v(K_t) = -v(t)$ for any additive invariant $v$ of order $n$. Then the ribbon knot $K\#2K_1\#2K_t$ is an $n$-inverse of $K$ in $\mathcal{RG}_n$. ∥

It is now immediate that finitely many additive Vassiliev invariants cannot distinguish the ribbon knots from the non-ribbon knots with trivial Arf invariant. The set of additive knot cobordism invariants and the set of additive Vassiliev invariants of order $> 2$ are disjoint. In other words, the Arf invariant is essentially the only Vassiliev invariant in the set of additive knot cobordism invariants. Corollary 3.5 and Corollary 3.6 can also be restated using $n$-equivalence.

## 5 An upper bound for the dimension of $V_n/V_{n-1}$

The evaluation of any primitive Vassiliev invariant on any split diagram is zero ([NS]) and any primitive invariant of order $n$, modulo the invariants of order $< n-1$, is determined by its values on chord diagrams of order $n$. Therefore, as a consequence of Theorem 6 in [B-N] and Theorem 1.2, any primitive invariant of order $n > 2$, modulo the invariants of order $< n$, is determined by its values on complete $n$-gons in $O_n$ and diagrams of order 2. Let $k$ be a field of characteristic not equal to two. Let $V'_n$ be the space of primitive Vassiliev invariants of order $\leq n$ over $k$ and let $V_n$ be the space of Vassiliev invariants of order $\leq n$ over $k$. Then an upper bound for the dimension of the quotient space $V'_n/V'_{n-1}$ may be obtained by counting the number of complete $n$-gons which generate $O_n$. In this section, we shall compute such an upper bound and then derive from it an upper bound for the dimension of the quotient space $V_n/V_{n-1}$. Unless otherwise stated, we assume that $n$ is an integer greater than 2.

Let $Y_n$ be the set of all $n$-cycles in the symmetric group $S_n$. Then $f : S_n \to O_n$ factors through the map $S_n \to Y_n$ which assigns to each $\sigma \in S_n$ the $n$-cycle $(\sigma(1) \ldots \sigma(n))$. Call the quotient map $f' : Y_n \to O_n$. By the antisymmetry relation, $f'(\sigma^{-1}) = (-1)^n f'(\sigma)$. Let $\tilde{O}_n$ be the quotient set of $O_n$ under the antisymmetry relations and let $X_n$ be the set of equivalence classes of $Y_n$ under the relation that $\sigma \sim \sigma^{-1}$. Then the map $f' : Y_n \to O_n$ induces a surjective map $X_n \to \tilde{O}_n$. Note that the cardinality of $\tilde{O}_n$ is an upper bound for $\dim(V'_n/V'_{n-1})$. We shall first compute this upper bound.

Let $G$ be the cyclic group generated by the cycle $t = (1 \ldots n) \in S_n$. Then the action of $G$ on $Y_n$ by conjugation induces a well-defined group action on $X_n$. It also induces through $f'$ a group action on $O_n$ by rotating the complete $n$-gons by multiples of $\frac{2\pi}{n}$. Let $\tilde{X}_n$ be the set of all orbits of $X_n$ under the action of $G$. Then the map $X_n \to \tilde{O}_n$ factors through the quotient map $X_n \to \tilde{X}_n$ to give a map $\tilde{f} : \tilde{X}_n \to \tilde{O}_n$ which is indeed a bijection.

The cardinality of $\tilde{O}_n$ is equal to the cardinality of $\tilde{X}_n$ which can be computed from Burnside's Lemma. Let $X(m)$ be the set of elements of $X_n$ fixed by $t^m$ and let $|S|$ denote the cardinality of a set $S$, then

$$|\tilde{X}_n| = \frac{1}{|G|} \sum_{0 \le m < n} |X(m)|$$

The trivial element $1 \in G$ fixes all elements in $X_n$, so $|X(0)| = |X_n| = \frac{(n-1)!}{2}$. It is easy to check that for any $0 < m < n$, $X(m) = X(d)$ where $d$ is the greatest common divisor of $m$ and $n$. Therefore, we only need to compute $|X(d)|$ for divisors $d$ of $n$.

**Lemma 5.1** *Let $0 < d < n$ and let $d$ divide $n$. Then*

$$|X(d)| = \begin{cases} \frac{1}{2}(d-1)! \left(\frac{n}{d}\right)^{d-1} \varphi(\frac{n}{d}) & \text{when } d \ne \frac{n}{2} \\ \frac{1}{2}(d-1)! \left(\frac{n}{d}\right)^{d-1} [\varphi(\frac{n}{d}) + d] & \text{when } d = \frac{n}{2} \end{cases}$$

*where $\varphi(k)$ is the number of positive integers not exceeding $n$ which are relatively prime to $n$.*

**Proof.** Let $\sigma = (a_1, \ldots, a_{n-1}, n) \in Y_n$ belong to an equivalence class in $X_n$ fixed by $t^d$. Then $t^d \cdot \sigma = \sigma$ or $\sigma^{-1}$.

Case 1. Suppose $t^d \cdot \sigma = \sigma$. Then $\sigma \cdot t^d = t^d$. Since $t^d$ is a product of $d$ disjoint cycles of length $\frac{n}{d}$, $\sigma$ acts on $t^d$ by permuting the $d$ cycles. Let $\Omega_j = \{j, d+j, \ldots, (\frac{n}{d}-1)d+j\}$ for $0 < j \leq d$. Then one can show that

**(O1)** $a_1, \ldots, a_d$ are all in distinct $\Omega_j$'s.

**(O2)** $a_k \equiv a_{d+k} + n - a_d \pmod{n}$ for $1 \leq k < n - d$.

**(O3)** The greatest common divisor of $a_d$ and $n$ is $d$.

Furthermore, the conditions (O1)-(O3) are sufficient for an $n$-cycle to commute with $t^d$. We have $(d-1)!$ ways of matching $\Omega_1, \ldots, \Omega_{d-1}$ with $a_1, \ldots, a_{d-1}$. Each $a_j$ ($1 \leq j < d$) has $|\Omega_j| = \frac{n}{d}$ choices while $a_d$ has $\varphi(\frac{n}{d})$ choices. Therefore, there are $(d-1)! \left(\frac{n}{d}\right)^{d-1} \varphi(\frac{n}{d})$ $n$-cycles which commute with $t^d$. Since the relation $\sim$ pairs up $\sigma$ and $\sigma^{-1}$, half of this number of elements of $X(d)$ satisfy Case 1.

Case 2. Suppose $t^d \cdot \sigma = \sigma^{-1}$. Then $t^{2d} \cdot \sigma = \sigma$. If the greatest common divisor of $2d$ and $n$ is $d$, then $t^d \cdot \sigma = \sigma$ contradicting that $t^d \cdot \sigma = \sigma^{-1}$. If $2d$ divides $n$, say $n = 2dl$. Let $a_k = (2l-1)d$. By comparing the two $n$-cycles $t^d \cdot \sigma$ and $\sigma^{-1}$, we get $a_k = d$ and $n = 2d$. The integer $k$ cannot be even. Moreover, the sets $\{a_{k+j}, a_{n-j}\}$ for $1 < j < n-k$ and $\{a_{k-j}, a_j\}$ for $1 \leq j < k$ form the classes $\Omega_1, \ldots, \Omega_d$. So for any odd $k$ with $1 \leq k \leq n$, if $a_1, \ldots, a_{\frac{k-1}{2}}, a_{k+1}, \ldots, a_{k+\frac{n-k-1}{2}}$ are all in distinct classes $\Omega_1, \ldots, \Omega_{d-1}$, then they determine an element of $X(d)$ in Case 2. There are $d$ choices for $k$ and given such an odd integer $k$, there are $(d-1)!2^{d-1}$ choices for such sequence of numbers. Again as in Case 1, the sequences occur in pairs. Thus, there are $\frac{1}{2}(d!2^{d-1})$ elements of $X(d)$ in Case 2.

In summary, if $n \neq 2d$, then all elements of $X(d)$ satisfy Case 1. If $n = 2d$, then there are $\frac{(d-1)!}{2} \left(\frac{n}{d}\right)^{d-1} \varphi(\frac{n}{d})$ elements of $X(d)$ satisfying Case 1 and $\frac{1}{2}(d!2^{d-1})$ elements of $X(d)$ satisfying Case 2. Thus, the lemma is proved. ∥

**Theorem 5.2** *Let $n$ be an integer greater than 2.*
*When $n$ is odd, $\dim(V'_n/V'_{n-1})$ is bounded above by*

$$\frac{1}{2n^2} \sum_{d|n} d! \left(\frac{n}{d}\right)^d \varphi(\frac{n}{d})^2.$$

*When $n$ is even, $\dim(V'_n/V'_{n-1})$ is bounded above by*

$$\frac{1}{2n^2} \sum_{d|n} d! \left(\frac{n}{d}\right)^d \varphi(\frac{n}{d})^2 + \frac{1}{4n} \left(\frac{n}{2}\right)! 2^{\frac{n}{2}}.$$

**Proof.** We have discussed that $\dim(V'_n/V'_{n-1})$ is bounded above by the cardinality of the quotient set $\tilde{X}_n$. Following the earlier discussions, we have

$$|\tilde{X}_n| = \frac{1}{n}\left(\sum_{\substack{d|n \\ d \neq n}} \varphi(d)|X(d)| + \frac{(n-1)!}{2}\right)$$

Then the theorem follows if we replace $|X(d)|$ by the formulas in Lemma 5.1.
∥

In Theorem 5.2, the term for $d = n$ is $\frac{n!}{2n^2}$ which is asymptotically equal to $\frac{1}{2}(n-2)!$. Observe that the term for $d = n$ dominates in both expressions in Theorem 5.2 for large $n$, so the upper bound is asymptotically equal to $\frac{1}{2}(n-2)!$. Indeed, after analysing the expressions, we obtain that the bound itself is bounded by $(n-2)!$ for $n > 2$. They give bounds better than the ones given in [CD1].

**Corollary 5.3** *Let $V_n$ be the space of all Vassiliev invariants of order at most $n$ over $k$. Then for any integer $n \geq 6$, $\dim(V_n/V_{n-1})$ is bounded above by $\frac{(n-2)!}{2}$.*

**Proof.** Let $p_n$ be the dimension of the quotient space $V'_n/V'_{n-1}$ and let $d_n$ be the dimension of the quotient space $V_n/V_{n-1}$. Then $p_1 = d_1 = 0$. As a $k$-algebra, the space $V_n$ is generated by the primitive invariants of order at most $n$ ([B-N],[G]). Therefore, we have

$$d_n = p_n + \sum_{\substack{k_1+\ldots+k_r=n \\ 2\leq k_1 \leq \ldots \leq k_r \leq n}} p_{k_1} \cdots p_{k_r} \qquad (2)$$

Theorem 5.2 already provides an upper bound for $p_n$'s. We shall first estimate those bounds. Observe that the function $\varphi(k)$ is bounded by $k$. In the summation over $d|n$, the term for $d = n$ is $\frac{(n-1)!}{2n}$. For $n \geq 12$ and $d \neq n$, $|X(d)|$ is bounded by $\frac{(n-3)!}{2}$. The number of divisors of $n$ that are less than $n$ is at most $\frac{n}{2}$. Combining these observations, we obtain that for $n \geq 12$,

$$\begin{aligned} p_n &\leq \frac{(n-1)!}{2n} + \frac{(n-3)!}{8} \\ &= \frac{(n-2)!}{2} - (1-\frac{2}{n})\frac{(n-3)!}{2} + \frac{(n-3)!}{8} \\ &= \frac{(n-2)!}{2} + (\frac{1}{n} - \frac{3}{8})(n-3)! \end{aligned} \qquad (3)$$

If we check the cases for small $n$, then the inequality (3) holds for $n \geq 6$ and $p_n \leq (n-2)!$ for $n \geq 2$. Thus, the summation over partitions of $n$ in (2) is bounded by

$$S = \sum_{r=2}^{[\frac{n}{2}]} \sum_{\substack{n_1+\ldots+n_r=n-2r \\ 0 \leq n_1 \leq \ldots \leq n_r}} n_1! \cdots n_r!$$

In each partition $(n_1, \ldots, n_r)$ of $n - 2r$ with $0 \leq n_1 \leq \ldots \leq n_r$, one must have $0 \leq n_j \leq \frac{n-2r}{r-j+1}$ for $1 \leq j \leq r$, and $n_1! \cdots n_r! \leq (n-2r)!$. Then for $n \geq 6$,

$$\begin{aligned}
S &\leq \sum_{\substack{n_1+n_2=n-4 \\ 0 \leq n_1 \leq n_2}} n_1! n_2! + \sum_{2 < r \leq [\frac{n}{2}]} (\frac{n-2r}{r}+1) \ldots (\frac{n-2r}{2}+1)(n-2r)! \\
&\leq (n-4)! + (\frac{n-4}{2})(n-5)! + \sum_{2 < r \leq \frac{n}{2}} \frac{(n-r)!}{r!(n-2r+1)} \\
&\leq \frac{3}{2}(n-4)! + \sum_{2 < r \leq \frac{n}{2}} \frac{2}{r!}(n-r-1)! \\
&\leq \frac{3}{2}(n-4)! + \{\frac{1}{3}(n-4)! + \frac{1}{12}(\frac{n}{2}-3)(n-5)!\} \\
&\leq 2(n-4)! \quad\quad\quad\quad\quad\quad\quad\quad\quad\quad\quad\quad\quad\quad\quad\quad\quad\quad\quad (4)
\end{aligned}$$

The result follows by applying the inequalities (3) and (4) to the equation (2). ∥

The first few numbers $(n > 2)$ computed from the formulas given in Theorem 5.2 are $1, 2, 4, 14, 54, 332, 2246$. The actual dimensions for $2 < n \leq 9$ computed by Bar-Natan [B-N] are $1, 2, 3, 5, 8, 12, 18$. Thus, our estimates are still quite rough. One may improve the upper bounds in Theorem 5.2 by finding more linear relations between the complete $n$-gons. For example, Figure 18 shows a relation for $n \geq 5$. In general, it is a hard problem to find all the relations and to get a linearly independent set which generates all diagrams.

Department of Mathematics
Columbia University
New York, NY10027